\documentclass[12pt]{article}
\usepackage{amsfonts}

\usepackage{amsmath,amssymb,amsthm}
\usepackage{multicol}
\usepackage{color}
\usepackage{hyperref}
\usepackage{graphicx}
\usepackage[english]{babel}

\newtheorem{Def}{Definition}[section]

\newcommand{\gf}{\mathfrak{g}}

\newcommand{\af}{\mathfrak{a}}

\newcommand{\sfr}{\mathfrak{s}}

\newcommand{\gfh}{\hat{\mathfrak{g}}}
\newcommand{\afh}{\hat{\mathfrak{a}}}

\begin{document}
\title{On affine extension of splint root systems}

\author{V.D.~Lyakhovsky$^1$, A.A.~Nazarov$^{1,2}$ \\
  {\small $^1$ Department of High-energy and elementary particle physics, SPb State University}\\
  {\small 198904, Saint-Petersburg, Russia,}
  {\small e-mail: lyakh1507@nm.ru}\\
  {\small$^{2}$ Chebyshev Laboratory,}
  {\small Department of Mathematics and Mechanics, SPb State University}\\
  {\small 199178, Saint-Petersburg, Russia}
  {\small email: antonnaz@gmail.com}}
\date{}
\maketitle

\begin{abstract}
Splint of root system of simple Lie algebra appears naturally in
the study of (regular) embeddings of reductive subalgebras. It can
be used to derive branching rules. Application of
splint properties drastically simplifies calculations of
branching coefficients. We study affine extension of splint root system of simple Lie algebra and obtain relations on theta and branching functions.
\end{abstract}

\section{Introduction}
\label{sec:introduction}

The term {\it splint} was introduced by D. Richter in \cite{richter2008splints} where the
classification of splints for simple Lie algebras was obtained. The fan $\Gamma \subset \Delta$ was
introduced in \cite{lyakhovsky1996rra} as a subset of root system describing recurrent properties of
branching coefficients for maximal embeddings. Injection fan is an efficient tool to study branching
rules. Later this construction was generalized to non-maximal embeddings and affine Lie algebras in
\cite{2010arXiv1007.0318L, ilyin812pbc}. In paper \cite{2011arXiv1111.6787L} we have shown that the
existence of a splint for a root system of a simple Lie algebra leads to simplifications in
reduction procedures of a Lie algebra module to modules of a subalgebra. This effect is based on the
injection fan and singular element properties of Lie algebra modules.

In the present note we discuss possible applications of splint in
a root system of simple Lie algebra related to representation
theory of affine Lie algebras. We discuss the structure of
injection fan for affine Lie algebras and show that it admits a
decomposition similar to that used in \cite{2011arXiv1111.6787L}
for simple Lie algebras. Such a decomposition leads to equations
for theta-functions. We study graded branching of affine Lie
algebra modules reduced to a finite-dimensional subalgebra and
discuss consequences of splint in this case.

\section{Splints and affine Lie algebras}
\label{sec:definitions}

Consider simple Lie algebra $\gf$ with a root system $\Delta$. Let
$\af_{1}\subset \gf$ be its reductive subalgebra of the same rank, such that
$\Delta_{\af}\equiv\Delta_{1}\subset\Delta$ and $Q_{\af}\subset Q$, where $Q$ is the root lattice. Irreducible highest-weight modules of
$\gf$ and $\af$ are denoted by $L^{\mu}$ and $L^{\nu}_{\af}$
correspondingly. Weyl character formula for irreducible modules is
$\mathrm{ch}L^{\mu}=\frac{\Psi^{(\mu)}}{\prod_{\alpha\in\Delta}(1-e^{-\alpha})}$,
where $\Psi^{(\mu)}=\sum_{w\in W}\epsilon(w) e^{w(\mu+\rho)-\rho}$
is a singular element of the module and $W$ -- the Weyl group of
$\gf$. Formal character of irreducible module admits a
decomposition
\begin{equation}
  \label{eq:4}
  \mathrm{ch} L^{\mu}=\sum_{\nu\in P_{\af}} b^{\mu}_{\nu} \mathrm{ch} L^{\nu}_{\af},
\end{equation}
where $P$, $P_{\af}$ are weight lattices of $\gf$ and $\af$. We
want to study affine extension of this situation:
$\gf\subset\gfh,\; \af\subset\afh$, $\afh\subset\gfh$,
$\hat{\Delta}\subset\hat{\Delta}$ and
$\mathrm{ch}L^{\hat{\mu}}_{\gfh}=\sum_{\hat{\nu}}
b^{\hat{\mu}}_{\hat{\nu}} \mathrm{ch} L^{\hat{\nu}}_{\afh}$. For
weights of an affine Lie algebra $\gfh$ we have
$\hat{\mu}=(\mu,k,n)$, where $\mu$ is a weight of $\gf$, $k$ --
the level of the module and $n$ -- the grade of the weight
$\hat{\mu}$

\begin{Def}
  Embedding $\phi$ of a root system $\Delta_1$ into a root system $\Delta$ is a bijective map of
  roots of $\Delta_{1}$ to a (proper) subset of $\Delta$ that commutes with vector composition law
  in $\Delta_{1}$ and $\Delta$.
\begin{equation*}
\phi:\Delta_1 \longrightarrow \Delta, \quad \phi \circ (\alpha + \beta) =\phi \circ \alpha + \phi \circ \beta,\,\,\, \alpha,\beta \in \Delta_1
\end{equation*}
\end{Def}

Note that the image $Im(\phi)$ must not inherit the root system properties except the addition rules
equivalent to the addition rules in $\Delta_{1}$ (for pre-images). Two embeddings $\phi_1$ and
$\phi_2$ can splinter $\Delta$ when the latter can be presented as a disjoint union of images
$Im(\phi_1)$ and $Im(\phi_2)$.

$\phi$ induces an injection of formal algebras $:{\mathcal{E}}_0
\hookrightarrow \mathcal{E}$ and for the image ${\mathcal{E}}%
_i=Im_{\phi}\left( {\mathcal{E}}_0\right)$ one can consider its inverse $%
\phi^{-1}:{\mathcal{E}}_i \longrightarrow {\mathcal{E}}_0$.

\begin{Def}
  A root system $\Delta $ ''splinters'' as $(\Delta _{1},\Delta _{2})$ if there are two embeddings
  $\phi _{1}:\Delta _{1}\hookrightarrow \Delta $ and $%
  \phi _{2}:\Delta _{2}\hookrightarrow \Delta $ where (a) $\Delta $ is the disjoint union of the
  images of $\phi _{1}$ and $\phi _{2}$ and (b) neither the rank of $\Delta _{1}$ nor the rank of
  $\Delta _{2}$ exceeds the rank of $%
  \Delta $.
\end{Def}

It is equivalent to say that $(\Delta_1,\Delta_2)$ is a "splint'' of $\Delta$ and we shall denote
this by $\Delta \approx (\Delta_1,\Delta_2)$. Each component $\Delta_1$ and $\Delta_2$ is a "stem''
of the splint.

We consider the case when one of the stems $\Delta _{1}=\Delta _{\frak{a}}$ is a root subsystem. As shown in paper \cite{2011arXiv1111.6787L} the second stem $\Delta _{\frak{s}}:=\Delta
_{2}=\Delta \setminus \Delta _{\frak{a}}$ can be translated into a product
$\prod_{\beta \in \Delta _{\frak{s}}^{+}}\left( 1-e^{-\beta }\right)=-\sum_{\gamma \in P}s(\gamma )e^{-\gamma }\quad   \label{splint product}$
and it defines an injection fan $\Gamma _{\frak{a}%
\hookrightarrow \frak{g}}$ \cite{lyakhovsky1996rra,ilyin812pbc,2010arXiv1007.0318L}.

Since the singular element of $L^{\mu}$ can be
written as $\Psi _{\frak{g}}^{\left( \mu \right) }=e^{-\rho }
\sum_{w\in W_{\frak{a}}}\epsilon \left( w \right) w\circ \left(
e^{\rho _{\frak{a}}}\Psi ^{\widetilde{\mu }+\rho
_{\frak{s}}}\right)$  for branching coefficients we get the
identity \cite{2011arXiv1111.6787L}:
\begin{equation}
  \label{eq:9}
b_{\left( \mu -\phi \left( \widetilde{\mu }-\widetilde{\nu }\right) \right)
}^{(\mu )}=M_{\left( \frak{s}\right) \widetilde{\nu }}^{\widetilde{\mu }}  
\end{equation}
Here the highest weight $\widetilde{\mu }$ is totally defined by the weight $\mu $, they have the
same Dynkin numbers: $\mu =\sum m_{k}\omega _{k}\qquad \Longrightarrow \quad \widetilde{\mu }=\sum
m_{k}\omega_{(\sfr) k} . \label{new h weight}$ So branching coefficients coincide with weight
multiplicities of $\sfr$-modules.

Now we consider affine extension of this setup, $\afh\subset\gfh$. Since $\mathrm{rank}\gf\leq
\mathrm{rank} \af+\mathrm{rank}\sfr$ for Weyl denominators we get
 \[
\prod_{\alpha\in\hat{\Delta}^{+}_{1}}(1-e^{-\alpha})^{\mathrm{mult}(\alpha)}\prod_{\beta\in\hat{\Delta}^{+}_{2}}(1-e^{\phi\circ \beta})^{\mathrm{mult}(\beta)}=\prod_{\gamma\in\hat{\Delta}^{+}}(1-e^{-\gamma})^{\mathrm{mult}(\gamma)}\prod_{n=0}^{\infty}(1-e^{-n\delta})^{\mathrm{rank}\af+\mathrm{rank}\sfr-\mathrm{rank}\gf}
\]

Using a specialization
\cite{kac1988modular,kac1984infinite,kac1990idl} and the
definition of Dedekind eta-function
$\eta(\tau)=q^{1/24}\prod_{n=1}^{\infty}(1-q^{n})$, where
$q=e^{2\pi i \tau}$ we can rewrite this identity as the relation
imposed on theta-functions $\Theta^{(\gfh)}_{\widehat{\lambda}=(\lambda,k,0)}(\tau,z)=\sum_{\xi\in Q_{\gf}+\frac{\lambda}{k}}e^{2\pi i k \left(\frac{1}{2} (\xi,\xi) \tau + (\xi,z)\right)}$:
\begin{multline}
  \label{eq:5}
  \eta(\tau)^{\mathrm{dim}(\af)}\prod_{\alpha\in\Delta_{1}^{+}}\frac{\Theta^{(\hat A_{1})}_{\alpha}(\tau,z)}{\eta(\tau)} \eta(\tau)^{\mathrm{dim}(\sfr)}\prod_{\beta\in\Delta_{2}^{+}}\frac{\Theta^{(\hat A_{1})}_{\phi\circ \beta}(\tau,z))}{\eta(\tau)}=\\
\eta(\tau)^{\mathrm{rank}(\af)+\mathrm{rank}(\sfr)-\mathrm{rank}(\gf)}
\eta(\tau)^{\mathrm{dim}(\gf)}\prod_{\alpha\in\Delta^{+}}\frac{\Theta^{(\hat A_{1})}_{\alpha}(\tau,z))}{\eta(\tau)}
\end{multline}
Here $z\in P_{\geq 0}\otimes \mathbb{C}$. Using Weyl denominator identity this relation can be
rewritten as a non-trivial relation connecting theta-functions of algebras $\gfh,\hat\sfr,\afh$:
\begin{equation}
  \label{eq:6}
  \left(\sum_{v\in W_{\af}}\epsilon(v) \Theta^{(\afh)}_{v\rho_{\af}}(\tau,z)\right)
  \cdot \left(\sum_{u\in W_{\sfr}}\epsilon(u) \Theta^{(\hat{\sfr})}_{\phi\circ(u\rho_{\sfr})}(\tau,z)\right)= 
  \left(\sum_{w\in W}\epsilon(w) \Theta^{(\gfh)}_{w\rho_{\gf}}(\tau,z)\right)
\end{equation}

Now consider the branching of $\gfh$-module to $\gf$-modules. For
formal characters we can write the following expression:
\begin{equation}
  \label{eq:8}
\mathrm{ch}L^{\hat{\mu}}_{\gfh}=\sum_{n=0}^{\infty}e^{-n\delta} \sum_{\nu\in P} b^{(\hat{\mu})}_{\nu}(n) \mathrm{ch} L^{\nu}_{\gf}.
\end{equation}
Rewriting this equation for weight multiplicities we get
$m^{(\hat{\mu})}_{\hat{\nu}=(\nu,k,n)}=\sum_{\xi\in P}
b^{(\hat{\mu})}_{\xi}(n) m^{(\xi)}_{\nu}$. We can introduce
branching functions similarly to the case of branching for affine
subalgebra \cite{kac1988modular,kac1990idl}:
$b^{(\hat{\mu})}_{\nu}(q)=\sum_{n=0}^{\infty}
b^{(\hat{\mu})}_{\nu}(n) q^{n}$.  These branching functions are connected to $q$-dimension of module $\mathrm{dim}_{q}L^{\hat \mu}_{\gfh}=\sum_{n=0}^{\infty}q^{n}\sum_{\nu\in P} b^{(\hat \mu)}_{\nu}(n) \mathrm{dim }L^{\nu}_{\gf}=\sum_{\nu\in P}b^{(\hat\mu)}_{\nu}(q) \mathrm{dim} L^{\nu}_{\gf}$. It is well-known that $q$-dimension is a modular function for some $\Gamma\subset SL_{2}(\mathbb{Z})$  \cite{gannon2006moonshine}, so branching functions $b^{(\hat \mu)}_{\nu}(q)$ have modular properties.

 For string functions of a module
$L^{\hat{\mu}}$ we have
\begin{equation}
  \label{eq:7}
   \sigma^{(\hat{\mu})}_{\nu}(q) = \sum_{\xi\in P} m^{(\xi)}_{\nu} b^{(\hat{\mu})}_{\xi}(q).
\end{equation}

Introduce an ordering of the set of weights $\xi$ as follows:
attribute to a weight $(\rho,\xi)$ its product $(\rho,\xi)$ with
the Weyl vector $\rho$. Then relation \eqref{eq:7} can be written
in the matrix form $\sigma(q)=M b(q)$ or as an inverse relation
$b(q)=M^{-1}\sigma(q)$. Here $\sigma(q)$ and $b(q)$ are infinite
columns of string and branching functions. Matrix $M$ contains
multiplicities of weights in $\gf$-modules similar to that of
Table 1 in paper \cite{2010arXiv1001}. The inverse matrix $M^{-1}$
encodes recurrent relations imposed weight multiplicities
\cite{il2010folded}.

Now consider the branching of $\gfh$-modules in $\af$-modules and
assume the existence of a splint
$\Delta^{+}_{\gf}=\Delta^{+}_{\af}\cup \phi(\Delta^{+}_{\sfr})$.
We decompose $\gf$-modules in equation \eqref{eq:8} into
$\af$-modules using property \eqref{eq:9}:
\begin{multline}
  \label{eq:10}
  \mathrm{ch}L^{\hat{\mu}}_{\gfh}=
\sum_{n=0}^{\infty}e^{-n\delta} \sum_{\nu\in P_{\af}} b^{(\hat{\mu})}_{(\gfh\downarrow\af)\nu}(n) \mathrm{ch} L^{\nu}_{\af}=
\sum_{n=0}^{\infty} e^{-n\delta} \sum_{\nu\in P} b^{(\hat{\mu})}_{(\gfh\downarrow\gf )\nu}(n) \sum_{\xi\in P_{\af}} b^{(\nu)}_{(\gf\downarrow \af) \xi}\mathrm{ch} L^{\xi}_{\af}=\\
=\sum_{n=0}^{\infty} e^{-n\delta} \sum_{\nu\in P} b^{(\hat{\mu})}_{(\gfh\downarrow\gf )\nu}(n) \sum_{\xi\in P_{\af}} M^{\widetilde{\nu}}_{  \widetilde{\nu}-\phi^{-1}( \nu-\xi )}\mathrm{ch} L^{\xi}_{\af}
\end{multline}
We see that the similar matrix relation holds for branching
functions $b_{(\gfh\downarrow\af)}(q)= M_{\sfr}\;
b_{(\gfh\downarrow\gf)}(q)$ and we can write
$\sigma(q)=M_{\af}\; b_{(\gfh\downarrow\af)}(q)$. So if we know
branching coefficients for the embedding $\gf\subset\gfh$ (for
example, see book \cite{kass1990ala}) we can easily obtain
branching functions for the embedding $\af\subset \gfh$.
\section*{Conclusion}
\label{sec:conclusion} We have demonstrated that splint in affine
Lie algebras leads to new relations between theta-functions and
branching functions for branching to finite-dimensional
subalgebras, which can be useful for computations. Further
question is to generalize this analysis to affine subalgebras and
to apply the results to branching in the study of CFT coset
models.

\section*{Acknowledgements}
\label{sec:acknowledgements}
Anton Nazarov thanks the Chebyshev Laboratory (Department of Mathematics and Mechanics,
Saint-Petersburg State University) for support under the grant 11.G34.31.0026
of the Government of Russian Federation.

\bibliography{bibliography}{}

\providecommand{\href}[2]{#2}\begingroup\raggedright\begin{thebibliography}{10}

\bibitem{richter2008splints}
D.~Richter, ``Splints of classical root systems,''
  \href{http://arxiv.org/abs/0807.0640}{{\tt arXiv:0807.0640}}.

\bibitem{lyakhovsky1996rra}
V.~Lyakhovsky, S.~Melnikov, {\em et al.}, ``{Recursion relations and branching
  rules for simple Lie algebras},'' {\em Journal of Physics A-Mathematical and
  General} {\bf 29} (1996) no.~5, 1075--1088,
  \href{http://arxiv.org/abs/q-alg/9505006}{{\tt q-alg/9505006}}.

\bibitem{2010arXiv1007.0318L}
V.~{Lyakhovsky} and A.~{Nazarov}, ``{Recursive algorithm and branching for
  nonmaximal embeddings},''
  \href{http://dx.doi.org/10.1088/1751-8113/44/7/075205}{{\em Journal of
  Physics A: Mathematical and Theoretical} {\bf 44} (2011) no.~7, 075205},
  \href{http://arxiv.org/abs/1007.0318}{{\tt arXiv:1007.0318 [math.RT]}}.

\bibitem{ilyin812pbc}
M.~Ilyin, P.~Kulish, and V.~Lyakhovsky, ``{On a property of branching
  coefficients for affine Lie algebras},'' {\em Algebra i Analiz} {\bf 21}
  (2009)  2, \href{http://arxiv.org/abs/0812.2124}{{\tt arXiv:0812.2124
  [math.RT]}}.

\bibitem{2011arXiv1111.6787L}
V.~{Laykhovsky} and A.~{Nazarov}, ``{Fan, splint and branching rules},''{\em
  ArXiv e-prints} (Nov., 2011)  , \href{http://arxiv.org/abs/1111.6787}{{\tt
  arXiv:1111.6787 [math.RT]}}.

\bibitem{kac1988modular}
V.~Kac and M.~Wakimoto, ``{Modular and conformal invariance constraints in
  representation theory of affine algebras},'' {\em Advances in mathematics(New
  York, NY. 1965)} {\bf 70} (1988) no.~2, 156--236.

\bibitem{kac1984infinite}
V.~Kac and D.~Peterson, ``{Infinite-dimensional Lie algebras, theta functions
  and modular forms},'' {\em Adv. in Math} {\bf 53} (1984) no.~2, 125--264.

\bibitem{kac1990idl}
V.~Kac, {\em {Infinite dimensional Lie algebras}}.
\newblock Cambridge University Press, 1990.

\bibitem{gannon2006moonshine}
T.~Gannon, {\em Moonshine beyond the Monster: The bridge connecting algebra,
  modular forms and physics}.
\newblock Cambridge Univ Pr, 2006.

\bibitem{2010arXiv1001}
M.~{Nesterenko}, J.~{Patera}, and A.~{Tereszkiewicz}, ``{Orthogonal polynomials
  of compact simple Lie groups},'' \href{http://arxiv.org/abs/1001.3683}{{\tt
  arXiv:1001.3683}}.

\bibitem{il2010folded}
M.~Il'in, P.~Kulish, and V.~Lyakhovsky, ``{Folded fans and string functions},''
  {\em Zapiski Nauchnykh Seminarov POMI} {\bf 374} (2010)  197--212.

\bibitem{kass1990ala}
S.~Kass, R.~Moody, J.~Patera, and R.~Slansky, {\em {Affine Lie algebras, weight
  multiplicities, and branching rules}}.
\newblock Sl, 1990.

\end{thebibliography}\endgroup
\bibliographystyle{utphys}

\end{document}